\documentclass[review]{elsarticle}

\usepackage{hyperref, mathtools}

\journal{Journal of \LaTeX\ Templates}

\DeclarePairedDelimiter\abs{\lvert}{\rvert}%

\hypersetup{
    colorlinks = true,
    linkcolor=blue
}

\newtheorem{conj}{Conjecture}
\newtheorem{theor}{Theorem}

\begin{document}

\begin{frontmatter}

\title{Some Open Problems Regarding the Number of Lines and Slopes in Arrangements that Determine Shapes}
\tnotetext[mytitlenote]{Built with Elsevier's article package, available on \href{http://www.ctan.org/tex-archive/macros/latex/contrib/elsarticle}{CTAN}.}

\author{Alexandros Haridis\corref{mycorrespondingauthor}}
\address{Department of Architecture, Massachusetts Institute of Technology}


\cortext[mycorrespondingauthor]{Corresponding author}
\ead{charidis@mit.edu}


\begin{abstract}

A set $L$ of straight lines and a set $P$ of points in the Euclidean plane define an arrangement $\mathcal{A} = (L, P)$ of construction lines and registration marks, if and only if: (1) any point in $P$ is a point of intersection of at least two lines in $L$, and (2) any two nonparallel lines in $L$ have a unique point of intersection in $P$. This paper discusses the following open problems regarding such arrangements. Suppose $k \geq 0$ number of points are given in the plane. How many construction lines $k$ points must determine? How many distinct slopes, or directions, are defined by construction lines that $k$ points determine? How many distinct sets of construction lines partition the plane, such that the lines meet at \emph{exactly} $k$ points? Empirical evidence is reported for small numbers of $k$, offering partial answers to the three problems. A conjecture is also stated for the first problem, on the number of construction lines, after examining a related problem about finite linear spaces from incidence geometry. This paper contributes to the body of work related to the \emph{mathematics of shapes} in the area of shape grammar theory.

\end{abstract}


\end{frontmatter}


\section{Overview of arrangements that determine shapes}

\subsection*{Shapes and their underlying arrangements}

In \cite{StinyPhD1975}, shapes are defined as finite sets of line segments in the two-dimensional Euclidean plane such that no two line segments in each set are adjacent co-linear or overlapping co-linear (such shapes are said to be in maximal element representation). For a brief overview of how shapes are defined, and for some background on notation and terminology, see \cite{HaridisGeomArxiv2020b}.

\begin{figure}[t!]
\centering
\includegraphics[width=\textwidth]{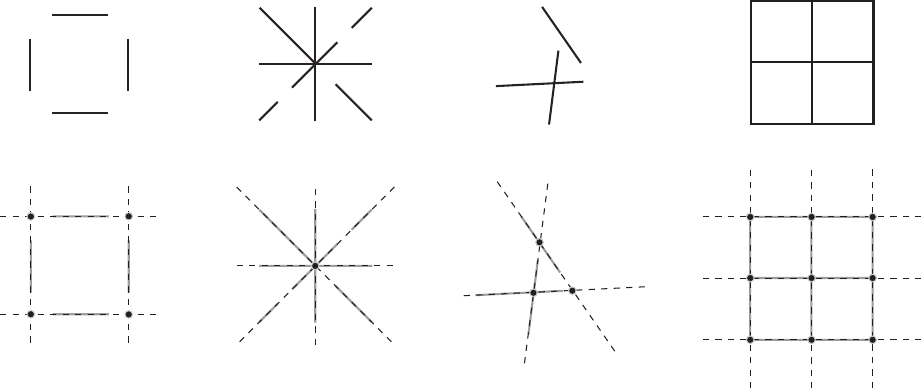}
\caption{Top: Examples of shapes. Bottom: The same shapes superimposed on their construction lines and registration marks.}
\label{Figure1}
\end{figure}

Given any shape, one can construct a special object out of the line segments of this shape called an \emph{arrangement} \cite{HaridisGeomArxiv2020b}. An arrangement associated with a shape $S$ is denoted by $\mathcal{A}_S$, and it consists of two things:

\begin{itemize}
\setlength{\parskip}{0pt}
\setlength{\itemsep}{0pt plus 1pt}
\item A set $L$ of straight lines called ``construction lines", and
\item A set $P$ of points called ``registration marks"
\end{itemize}

\noindent Construction lines are infinitely extending lines in the plane that contain one or more line segments of the shape (co-linear line segments separated by ``gaps" are contained in the same construction line). Registration marks are the points of intersection of pairs of construction lines. For example, Figure 1 shows some examples of shapes in the plane and their underlying arrangements. The arrangement $\mathcal{A}_S$ is defined as the pair $(L, P)$.

In this context, the lines in $L$ and the points in $P$ are defined in the ordinary way. A construction line is calculated from the line segments of the shape, using the familiar equation of lines. A registration mark is defined with specific coordinates and represents a solution to (at least) two equations of lines. Thus, both the shapes and their underlying arrangements are spatial geometric objects that can be executed/realized in the Euclidean plane.

If the sets $L$ and $P$ are empty, then a special object is defined called the \emph{empty arrangement}, denoted by $\mathcal{A}_{\emptyset}$, which is the arrangement associated with the \emph{empty shape}---a shape with no line segments.

It is common to use the term ``arrangement" to refer only to the set of lines $L$ (e.g. as in hyperplane arrangements). Here, the term ``arrangement" is used in place of the longer ``arrangement of construction lines and registration marks," and it refers to both the lines and the planes, that is, to both the sets $L$ and $P$. The term ``arrangements for shapes" is also used to mean the same thing.

Two or more shapes can have the same underlying arrangement. More generally, any arrangement of construction lines and registration marks determines a class of uncountably many equivalent shapes in the plane. For example, the shapes in Figure 2 are equivalent in the sense that they have the same underlying arrangement---the third one in Figure 1.

\begin{figure}[t!]
\centering
\includegraphics[width=\textwidth]{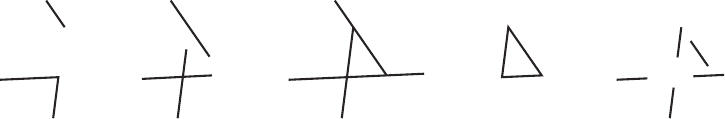}
\caption{Examples of shapes with the same underlying arrangement of construction lines and registration marks.}
\label{Figure2}
\end{figure}

\subsection*{Rules for recognizing arrangements of construction lines and registration marks}

It is a straightforward exercise to construct an arrangement from the line segments of a given shape \cite{HaridisGeomArxiv2020b}. In the opposite direction, given a set $L$ of straight lines and a set $P$ of points in the plane, how do we decide if the two sets form a valid arrangement of construction lines and registration marks? In other words, how do we decide if the two sets define an arrangement that contains shapes?

\begin{figure}[t!]
\centering
\includegraphics[width=\textwidth]{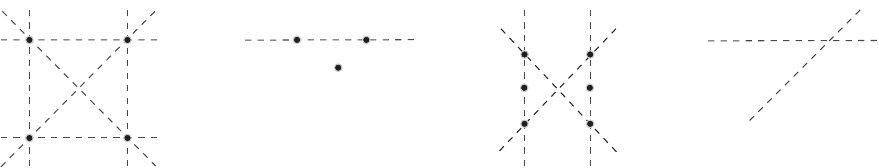}
\caption{Some sets of points and lines that do not define arrangements of construction lines and registration marks.}
\label{Figure3}
\end{figure}

A set $L$ of straight lines and a set $P$ of points in the Euclidean plane given in a standard coordinate system, define an arrangement $\mathcal{A} = (L, P)$ of construction lines and registration marks, if and only if:

\begin{enumerate}[(1)]
    \setlength{\parskip}{0pt}
    \setlength{\itemsep}{0pt plus 1pt}
    \item any point in $P$ is a point of intersection of at least two lines in $L$, and 
    \item any two nonparallel lines in $L$ have a unique point of intersection in $P$.
\end{enumerate}

The two rules, (1) and (2), do not assume that there are any lines at all. In this case, there cannot be points, and we get an \emph{empty arrangement}. If there are lines, the rules do not imply that there are points because the lines may be parallel. On the other hand, if there are points, the rules imply that there must be lines that meet at those points.

From a practical standpoint, the first rule guarantees there are no ``isolated" points (i.e. points not on any line) or points on just one line. The second rule guarantees all intersections of pairs of lines are in the set $P$. Figure 3 shows examples of sets of points and lines in the plane that do not define arrangements for shapes. Each pair of sets of points and lines fails to satisfy one or both rules.

\section{Line and slope problems}

The construction lines form a partition of the plane while the registration marks ``punctuate" the locations/points in the plane at which these lines meet. Suppose $k \geq 0$ number of points are given in the plane. Then, the following questions arise naturally: 

\begin{itemize}
\setlength{\parskip}{0pt}
\setlength{\itemsep}{0pt plus 1pt}
\item How many construction lines $k$ points must determine?
\item How many distinct slopes, or directions, are defined by construction lines that $k$ points determine?
\item How many distinct sets of construction lines partition the plane, such that the lines meet at \emph{exactly} $k$ points?
\end{itemize}

\noindent The first two problems are problems about bounds. The third problem is an enumeration problem. 

The term ``distinct" depends on the particular interpretation it is given in a certain context. Two arrangements, for example, can be ``projectively distinct" when there is no affine or, more generally, projective transformation that makes one an image of the other under that transformation. Two arrangements can also be ``combinatorially distinct" when there is no way of distinguishing the arrangements solely on the basis of information about incidences between their respective sets of points and lines. The two methods of comparing arrangements are described in more detail in \hyperref[appendixB]{Appendix B}.

To the extent of my knowledge, no results exist for the above three problems that specifically refer to arrangements of construction lines and registration marks satisfying the two rules given in the previous section. There are related problems in the literature of finite combinatorial geometry (e.g. \cite{BattenCombinatoricsFiniteGeom, Grunbaum2009, AignerZiegler2018})---one of these is discussed in Section 4---but these problems assume that arrangements of points and lines are not necessarily spatial geometric objects but more general abstract incidence systems (e.g. \cite{BattenCombinatoricsFiniteGeom}), satisfying incidence laws which are different from those that the arrangements we're interested in here satisfy (this distinction is also discussed in \cite{HaridisGeomArxiv2020b}). 

In this paper, empirical evidence is reported for small numbers of $k$ offering partial answers to the three problems above. In Section 4, a conjecture is stated for the first problem on the number of construction lines, after examining a related problem about finite linear spaces from incidence geometry.

\section{Empirical evidence}

When $k = 0, 1, 2$ or $3$ points are given in the plane, it is a straightforward exercise to find bounds on the number of construction lines the $k$ points must define and to describe and visualize each distinct set of construction lines that partitions the plane.

\paragraph{k = 0} If $k = 0$, an arrangement contains no points and therefore it must fall under (either) one of the following two cases: (i) the empty arrangement, with no points and no lines, (ii) an arrangement with only one construction line or with two or more parallel construction lines. In the second case, there is no upper bound on the number of construction lines defined although any arrangement must have \emph{finitely} many points and lines. In case (ii) a single slope is defined. The two cases are shown graphically in this figure.

\begin{figure}[h!]
\centering
\includegraphics{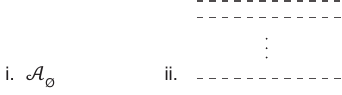}
\label{Figure4}
\end{figure}

\paragraph{k = 1} If $k = 1$ arrangements contain one point and a finite number of construction lines greater than or equal to two all of which intersect at that point. If $n$ is the number of construction lines ($2 \leq n < \infty$), then $n$ distinct slopes are defined by those lines.

\begin{figure}[h!]
\centering
\includegraphics{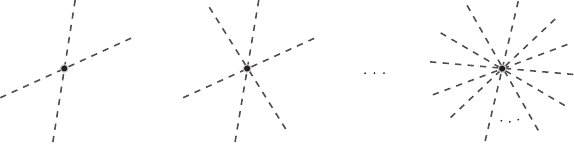}
\label{Figure5}
\end{figure}

\paragraph{k = 2} If $k = 2$, an arrangement contains two points and exactly three construction lines two of which are parallel and a third that intersects both once. Exactly two distinct slopes are defined.

\begin{figure}[h!]
\centering
\includegraphics{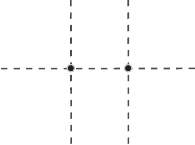}
\label{Figure6}
\end{figure}

\paragraph{k = 3} If $k = 3$, an arrangement contains three points and therefore it falls under (either) one of the following three cases: (i) arrangements with four construction lines, three of which are parallel and a fourth one that intersects each of them once, (ii) arrangements with three pair-wise intersecting construction lines forming a triangle, with an additional fourth construction line passing through one (and only one) of the points of the triangle that's parallel to one of the other three construction lines, and (iii) arrangements with three pair-wise intersecting construction lines forming a triangle. In case (i), exactly two distinct slopes are defined; in (ii) and (iii), exactly three distinct slopes are defined. The three cases, (i), (ii), and (iii), are shown graphically in this figure.

\begin{figure}[h!]
\centering
\includegraphics{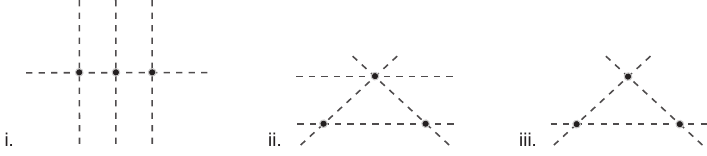}
\label{Figure7}
\end{figure}

\vspace{12pt}

As the number $k$ increases, it becomes significantly harder to describe arrangements verbally or visualize them one by one. Following an empirical experimentation by hand, a catalog has been produced which enumerates distinct arrangements of construction lines and registration marks in the plane for $k = 4, 5, 6$, and $7$. This catalog is shown in \hyperref[appendixA]{Appendix A}. Every arrangement in the catalog is marked by two signatures. The first signature is a triple of numbers, $(k\;\;n \;\; s)$, where $k = |P|$ the number of registration marks, $n = |L|$ the number of construction lines, and $s = |D|$ the number of distinct slopes/directions defined by the lines (see also \hyperref[appendixB]{Appendix B} ). The second signature encodes information about the number of lines meeting at each intersection point, and is given in this form $[2^{\# k2} 3^{\# k3} ...]$. The entries indicate ``the number of points in which 2 lines intersect, then the number of points in which 3 lines intersect, and so on". In this way, the sum of the superscripts is always equal to the total number of registration marks in an arrangement. Because every registration mark must be incident with at least two construction lines, all signatures must start with $[2^{\# k2} ...]$. We omit any entry of the form $c^0$. 

\section{Discussion}

\subsection*{Three common arrangement types}

There are three types of arrangements of construction lines and registration marks that are defined for any $k \geq 3$ number of registration marks:

\begin{itemize}
    \setlength{\parskip}{0pt}
    \setlength{\itemsep}{0pt plus 1pt}
    \item \emph{Near-pencil}: one line has $k - 1$ points and all others have two points, connecting the $k$th point with each of the $k - 1$ points. 
    \item \emph{Augmented near-pencil}: a near pencil arrangement, with one additional construction line passing through the $k$th point, parallel to the line containing the $k - 1$ points.
    \item \emph{Railtrack}: one line has $k$ points and all others have exactly one point, the point where they intersect the line with $k$ points.
\end{itemize}

The three types are shown graphically in Figure 4. They can be easily identified in the catalog in \hyperref[appendixA]{Appendix A}.

\begin{figure}[t!]
\centering
\includegraphics{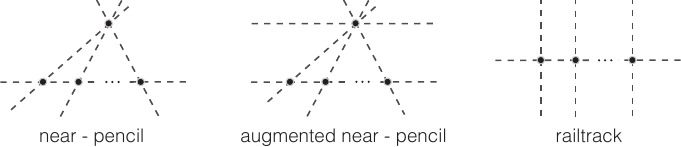}
\caption{Three types of arrangements that exist for any number of points greater than or equal to three.}
\label{Figure8}
\end{figure}

\subsection*{The line problem of de Bruijn and Erd{\H o}s}

Consider the following well-known theorem about finite linear spaces by de Bruijn and Erd{\H o}s (e.g. \cite{BattenCombinatoricsFiniteGeom, AignerZiegler2018}). 

\begin{theor}
\label{th1}
(de Bruijn-Erd{\H o}s, 1948)\;\;\;Let $P$ be a set of $k \geq 3$ points, and $L$ a set of $n > 1$ lines such that any two distinct points of $P$ are on exactly one line of $L$. Then

\begin{itemize}
    \setlength{\parskip}{0pt}
    \setlength{\itemsep}{0pt plus 1pt}
    \item $n \geq k,$
    \item if $n = k$, either one line has $k - 1$ points and all others have two points, or every line has $k + 1$ points and every point is on $k + 1$ lines, $k \geq 2$.
\end{itemize}
\end{theor}

In this theorem, it is assumed that ``lines" in the set $L$ are subsets of at least two points of $P$. And since any two distinct points in $P$ are on a unique line in $L$, the sets $P$ and $L$ form a linear space; for the rules that determine linear spaces see \cite{BattenCombinatoricsFiniteGeom}, and for a discussion of their relationship to arrangements of construction lines and registration marks, see \cite{HaridisGeomArxiv2020b}. 

\hyperref[th1]{Theorem 1} is not restricted to the geometric setting of points and lines in the Euclidean plane. It applies to any abstract system consisting of elements called ``points" and ``lines" which obey the desired ``incidence" rules (e.g. any two ``points" are on a ``line"). The meaning of ``points" and ``lines" depends on the interpretation given to these notions in the context of a problem or application.

Given $k \geq 3$ points in the plane not all on a line, the theorem provides a concrete lower bound on the number of lines passing through at least two points: \emph{$k$ points must determine at least $k$ different lines} \cite{AignerZiegler2018}. And in particular, according to the second bullet, when exactly $k$ lines are determined the sets $P$ and $L$ must form a near-pencil or a projective plane \cite{BattenCombinatoricsFiniteGeom}. 

This result by de Bruijn and Erd{\H o}s, is closely related to the first problem stated in Section 2: how many different construction lines $k$ registration marks must determine? However, the place where arrangements of construction lines and registration marks depart from the point-line incidence systems that \hyperref[th1]{Theorem 1} applies to, is in the assumption that ``any two points in $P$ are on a line in $L$". Which arrangements of construction lines and registration marks satisfy this assumption? Put in another way, which arrangements form a (finite) linear space?

An arrangement of construction lines and registration marks that forms a linear space is the near-pencil

\begin{figure}[h!]
\centering
\includegraphics{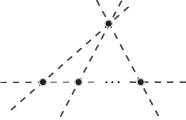}
\label{Figure9}
\end{figure}

\noindent A near-pencil can be defined for any $k \geq 3$ number of registration marks. A near-pencil with $k = 3$ points is a triangle---the smallest possible linear space. By \hyperref[th1]{Theorem 1}, a near-pencil arrangement with $k$ registration marks determines exactly $k$ construction lines (this is also easy to verify empirically, from the catalog in \hyperref[appendixA]{Appendix A}).

Empirical evidence also indicates that for $k = 5$ another linear space is defined, this one in particular

\begin{figure}[h!]
\centering
\includegraphics{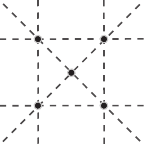}
\label{Figure10}
\end{figure}

\noindent This five-point arrangement determines $k + 1 = 6$ different lines.

I have not been able to identify other arrangements containing shapes that define linear spaces (i.e., in which any two distinct points lie on a line).

The other two common types of arrangements, apart from near-pencils, namely, augmented near-pencils and railtracks

\begin{figure}[h!]
\centering
\includegraphics{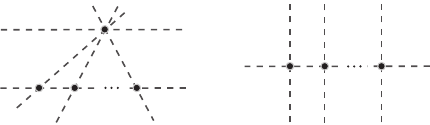}
\label{Figure11}
\end{figure}

\noindent also determine exactly $k + 1$ construction lines by definition. They can be defined for any $k \geq 3$, and both have exactly one construction line more than a near-pencil with the same number of points. However, neither of these two arrangement types defines a linear space: they contain lines incident with one point only. More generally, neither of these two types of arrangements defines a point-line geometry \cite{HaridisGeomArxiv2020b}.

It is reasonable to ask then, are there any arrangements containing shapes with $k$ registration marks such that they determine greater than $k + 1$ construction lines? Such arrangements, if they indeed exist, must necessarily determine \emph{more} lines connecting ``pairs of distinct points" than the near-pencils and they must do so without violating the two rules given in the definition of arrangements. We conjecture that such arrangements do not exist:

\begin{conj}
\label{conj1}
At most $k + 1$ construction lines are determined by $k \geq 3$ registration marks in the plane.
\end{conj}

Experiment seems to indicate that this is true, but a formal proof is overdue. A positive answer to \hyperref[conj1]{Conjecture 1} would provide an upper bound on the number of construction lines determined by $k \geq 3$ registration marks. No comparable conjecture or a concrete result is known for a lower bound on the number of construction lines.

One way to validate the above conjecture is the following. It is enough to show that an arrangement with $k > 5$ registration marks is a linear space, if and only if, it is a near-pencil. For $k = 3, 4, 5$, the near-pencils are the only linear spaces with the exception of $k = 5$ where there is an isolated case of a linear space given just above.

For the second problem in Section 2, if the above conjecture is true, $k \geq 3$ registration marks in the plane must determine \emph{at most} $k$ distinct slopes (the slopes determined by near-pencils), and \emph{at least} 2 distinct slopes. The lower bound follows from the rule that defines registration marks, according to which a registration mark is at the intersection of at least two construction lines.

Finally, for the third problem in Section 2, an algorithmic approach seems the most reasonable path toward enumerating the distinct arrangements of construction lines that partition the plane for given values of $k$.

The point-line arrangements studied in this paper contribute to the area of mathematics of shapes developed within shape grammar theory. Concurrent studies on the topic of finite topology and topological descriptions for shapes that represent designs can be found in \cite{HaridisEPB2020} and \cite{HaridisJMA2020}. A related study investigates how description rules encoding natural language descriptions of shape rules can be used as a method to provide both a constructive and an evocative description of a computation \cite{HaridisSpatialRulesArxiv2021}. 

\section*{Acknowledgement}
\noindent The content of this paper has been revised and subsumed in the following publication:

\vspace{0.2in}

\noindent Alexandros, Haridis. 2024. ``Arrangements containing shapes: mathematical features and their use in visual calculating." In SD Kotsopoulos (Ed.) \emph{Shape Computation: Fifty Years 1972-2022} (Springer Nature). Also in \emph{Mathematics and the Built Environment} book series.



\bibliography{sample}

\clearpage
\appendix
\section{Distinct arrangements of construction lines with $k = 4, 5, 6, 7$ registration marks}
\label{appendixA}

\begin{figure}[h!]
\includegraphics[width=4.4in]{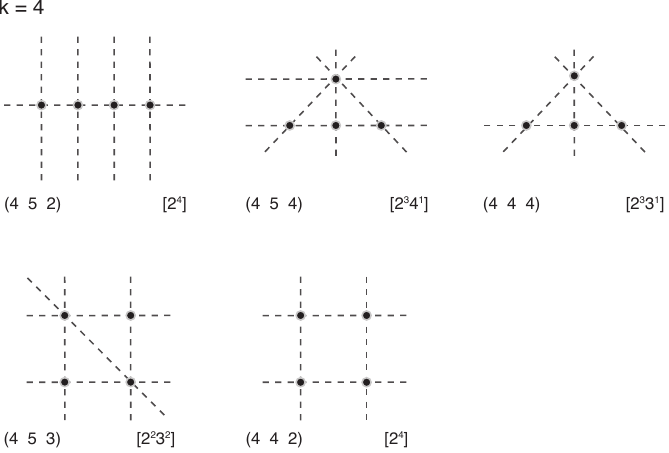}
\label{Figure12}
\end{figure}

\begin{figure}[h!]
\includegraphics[width=\textwidth]{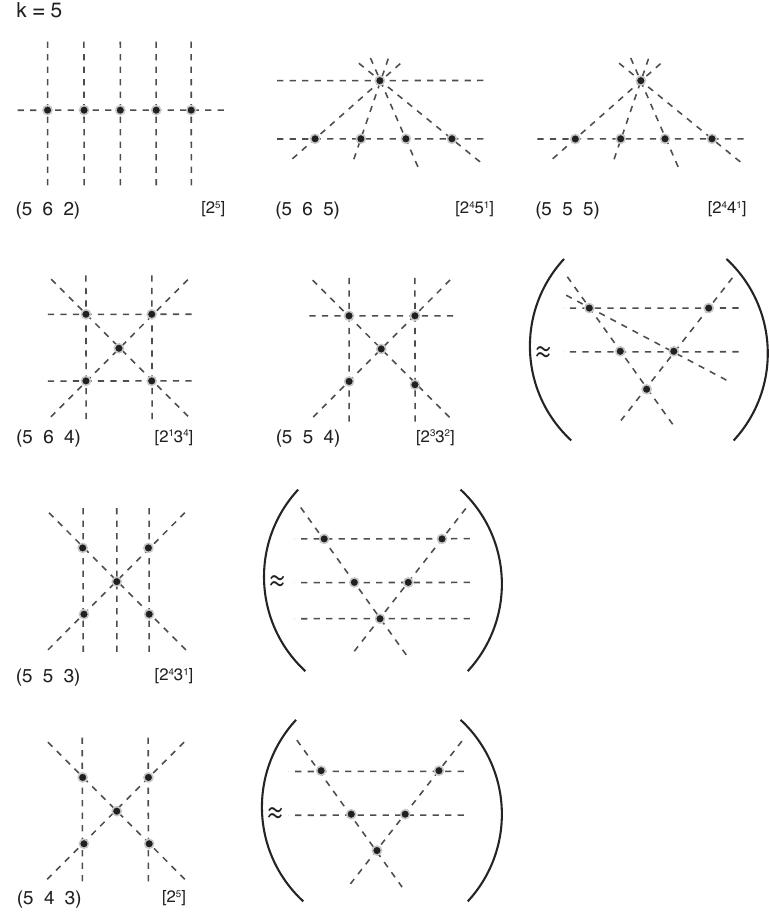}
\label{Figure13}
\end{figure}

\begin{figure}[h!]
\includegraphics[width=\textwidth]{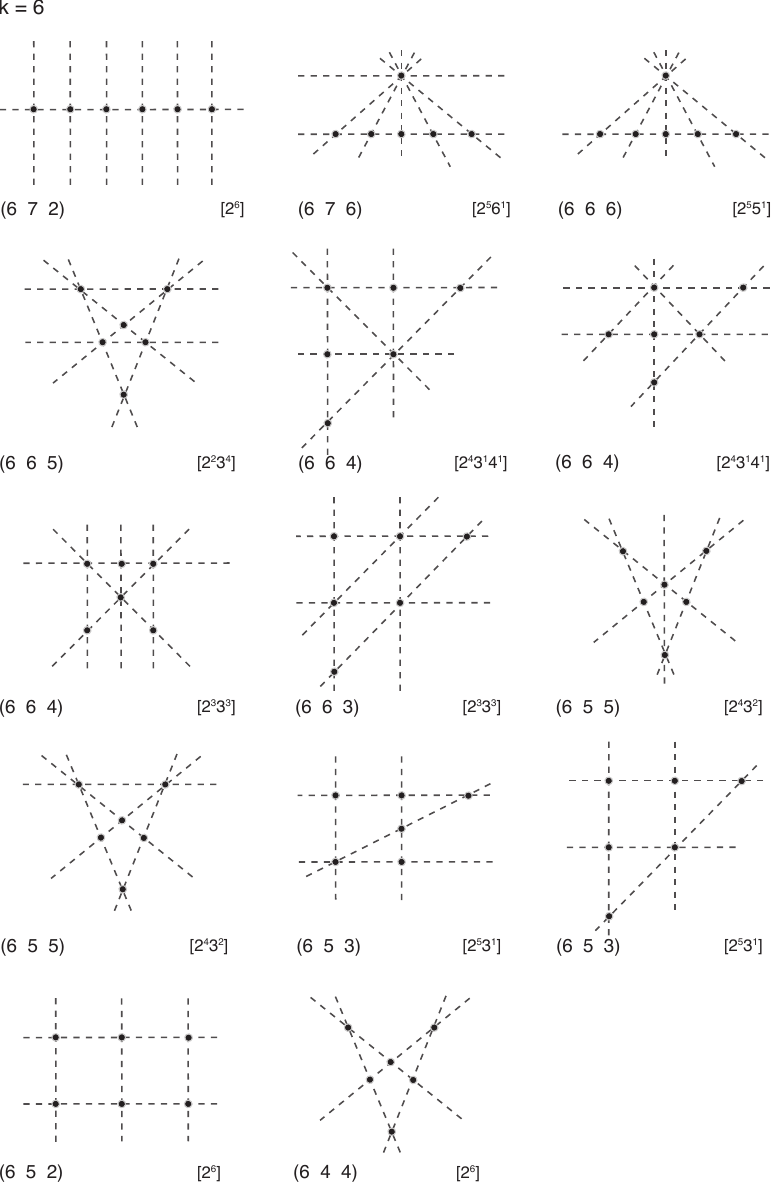}
\label{Figure14}
\end{figure}

\begin{figure}[h!]
\includegraphics[width=\textwidth]{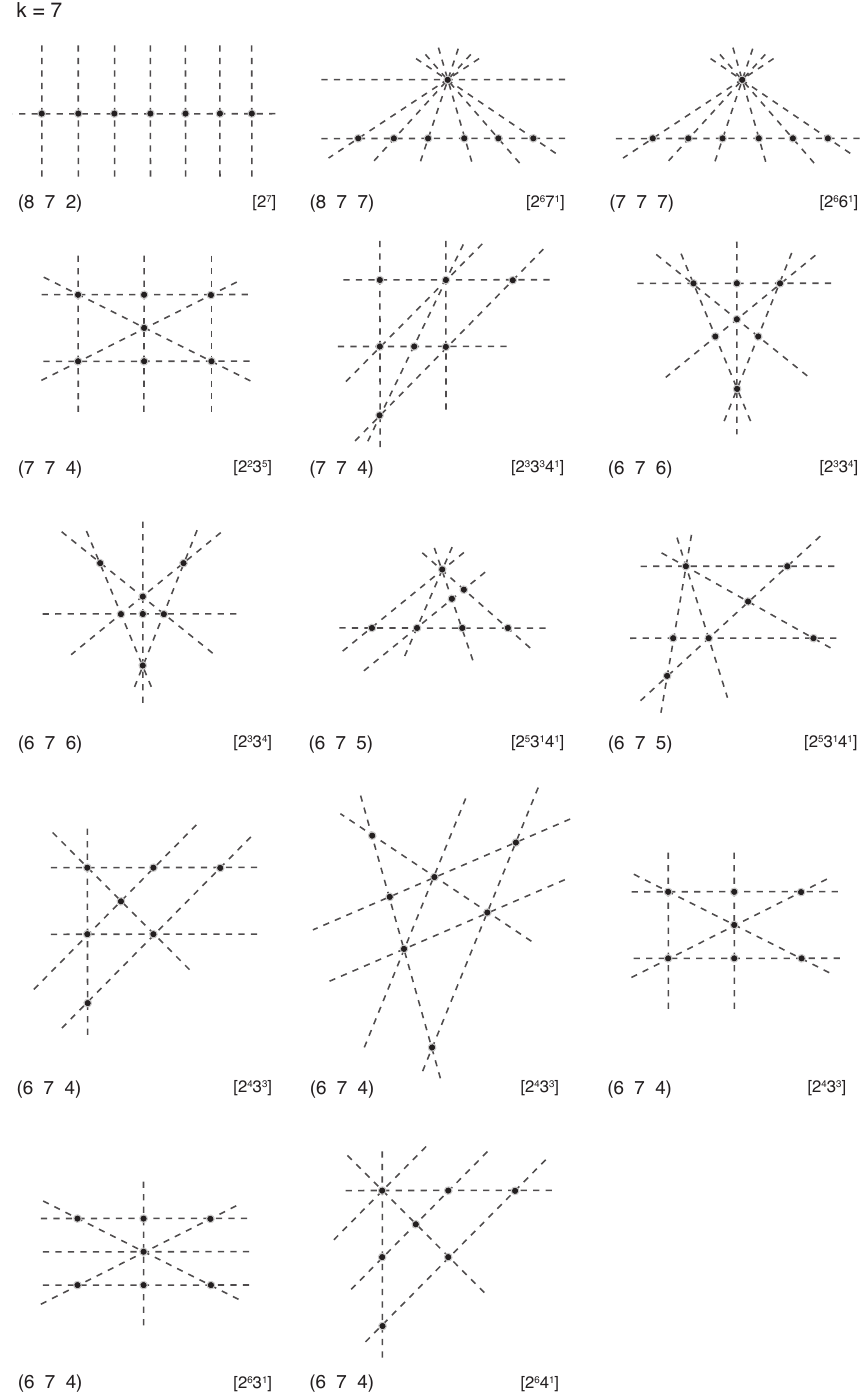}
\label{Figure15}
\end{figure}

\begin{figure}[h!]
\includegraphics[width=\textwidth]{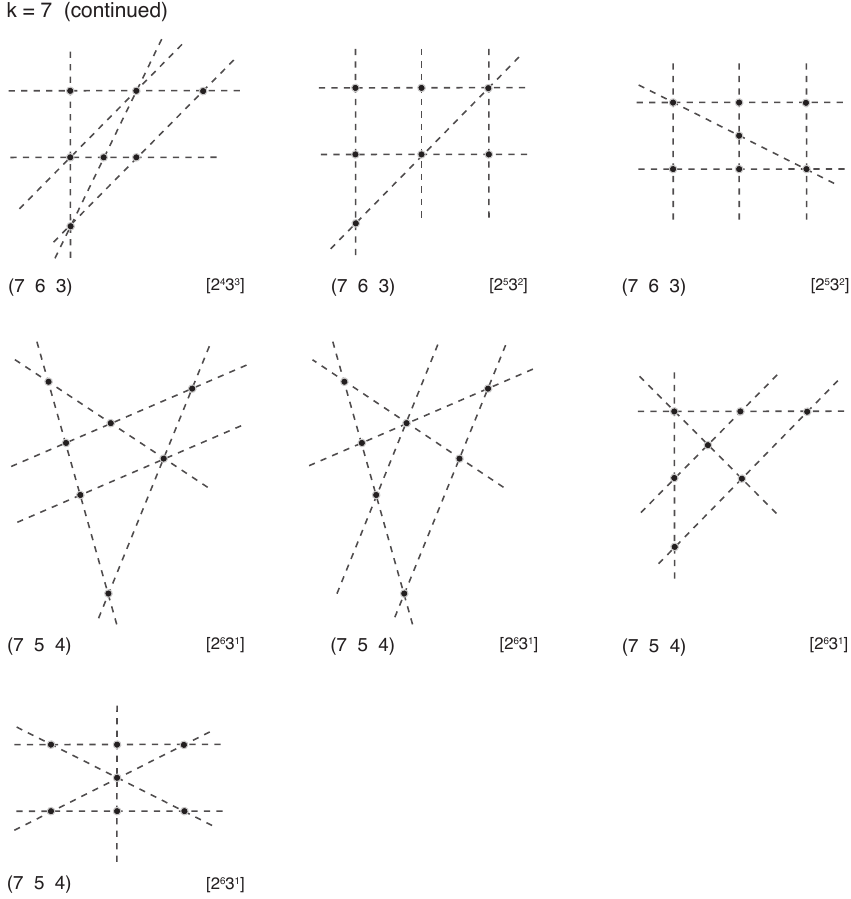}
\label{Figure16}
\end{figure}

\clearpage

\section{Comparing arrangements}
\label{appendixB}

\subsection{Projectively equivalent arrangements}

The physical realization (drawing) of an arrangement is in general not unique. The same arrangement of registration marks and construction lines can be drawn/represented in multiple ways (indeed, infinitely many). If the drawings of two arrangements differ only by a (non-singular) affine or projective transformation, the two arrangements are called \emph{projectively equivalent}. That is to say, their drawings represent the ``same" arrangement under a projective equivalence. For example, these three arrangements are equivalent because they are obtained as images of each other under a simple affine transformation.

\begin{figure}[h!]
\includegraphics[width=\textwidth]{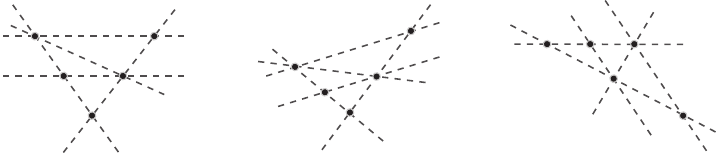}
\label{Figure17}
\end{figure}

Arrangements whose drawings cannot be mapped onto each other by an affine or projective transformation are called \emph{projectively distinct}. Visually, the drawings of such arrangements are expected to be quite dissimilar. However, even if two arrangements are projectively distinct in the aforementioned sense, they may still be equivalent in terms of the \emph{incidence relation} they determine between points and lines (for more information on incidence relations, see \cite{HaridisGeomArxiv2020b}).

\subsection{Combinatorially equivalent arrangements}

Projectively distinct arrangements can be distinguished from each other in terms of the following comparison measures.

\begin{enumerate}[(1)]
    \setlength{\parskip}{0pt}
    \setlength{\itemsep}{0pt plus 1pt}
    \item \emph{Number of registration marks, construction lines}, and \emph{distinct slopes}
    \item \emph{Line type}, \emph{slope type}, and  \emph{point type}
    \item \emph{Point-line degree type}
    \item \emph{Central degree and central signature}
\end{enumerate}

\vspace{12pt}

\paragraph{1. Number of registration marks, construction lines, and distinct slopes}

The simplest way to separate arrangements, is to compare their number of registration marks, construction lines, and distinct slopes (or directions).

Every construction line in the set $L$ has a slope by definition, denoted by $s(\mathcal{L})$, and two construction lines determine the same slope if and only if they are parallel. Thus, to any set $L$ of construction lines, there is a finite set

\vspace{12pt}
\centerline{$D = \{\;s(\mathcal{L})\;|\;\mathcal{L} \in L\;\}$}
\vspace{12pt}

\noindent of \emph{distinct slopes} of these lines. The number of distinct slopes in an arrangement is given by the size of the set $D$. Notice that $\abs{D} \leq \abs{L}$, for any arrangement.

Any given arrangement can be described as a triple of three simple numbers

\vspace{12pt}
\centerline{$(n\;\;k \;\; s)$}
\vspace{12pt}

\noindent where $n = |L|$, $k = |P|$, and $s = |D|$. 

The values for $n$, $k$, and $s$ do not determine arrangements uniquely---multiple arrangements can share the same values. However, they provide a first (crude) way of separating projectively distinct arrangements. For example, see \hyperref[appendixA]{Appendix A} where each arrangement is marked by such a triple of numbers.

\paragraph{2. Line type, slope type, and point type}

The \emph{line type} of an arrangement is given by a vector $(t_1, t_2, ...)$, where $t_u$, $u \geq 1$, is the number of lines which are incident with exactly $u$ points. The line type is also called the $t$-$vector$ of an arrangement. For example, the line type of all three arrangements in the figure in Section B.1 is $(0, 3, 2)$.

The \emph{slope type} of an arrangement is given by a vector $(s_1, s_2, ...)$, where $s_v$, $v \geq 1$, is the number of distinct slopes of $v$-point lines (i.e. lines passing through $v$ number of points). The slope type is also called the $s$-$vector$ of an arrangement. For example, the slope type of all three arrangements in the figure in Section B.1 is $(0, 2, 2)$.

The \emph{point type} of an arrangement is given by a vector $(p_2, p_3, ...)$ where $p_j$, $j \geq 2$, is the number of points incident with exactly $j$ lines, i.e. points of $j$-degree. The slope type is also called the $p$-$vector$ of an arrangement. For example, the point type of all three arrangements in the figure in Section B.1 is $(3, 2)$.

\vspace{12pt}

When the three vectors have unequal sizes, it is possible to use a simple technique of padding to enforce the same size for all three vectors---this may be of practical use if one wants to compare arrangements using vector-based metrics.

When two arrangements have the same line type, then they can be further distinguished in terms of their slope type. If their slope types are the same, then they can be distinguished in terms of their point type. When all three type vectors are the same, then the arrangements are equivalent on the basis of those vectors (they cannot be distinguished in terms of them). In this case, further means of diversification are required.

\paragraph{3. Point-line degree type}

The number of lines passing through a point $p$ is called the \emph{degree} of the point, and is denoted by $[p]$. The number of $u$-point lines that pass through $p$ is denoted by $[p]_u$, $u \geq 1$. If there are $k$ points in total in an arrangement, the \emph{point-line degree type} of a point $p$ is the vector $([p]_1, [p]_2, ..., [p]_k)$. For example, in the following figure the point-line degree type of point $q$ is $(0, 2, 1, 0, 0)$.

\begin{figure}[h!]
\centering
\includegraphics{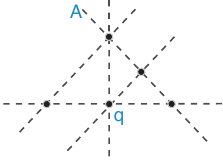}
\label{Figure18}
\end{figure}

The comparison measures in (2) and (3), can be employed to distinguish arrangements when their respective values for $n$, $k$, and $s$ are identical. 

If two arrangements are identical in terms of the measures in (2) and (3), then these arrangements have the same \emph{incidence relation}. Such arrangements are isomorphic or combinatorially equivalent.

An alternative way of identifying isomorphic arrangements is by labelling: two arrangements have the same incidences, provided their points and lines can be given such labels that a point and a line are incident in one of them if and only if they are incident in the other. This is illustrated in the following figure.

\begin{figure}[h!]
\centering
\includegraphics[width=\textwidth]{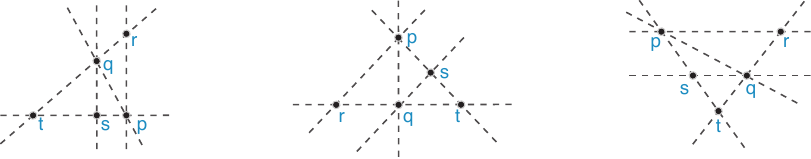}
\label{Figure19}
\end{figure}

\paragraph{4. Central degree and central signature}

Two isomorphic arrangements cannot be distinguished in terms of point-line incidences. However, isomorphic arrangements can be distinguished on the basis of a local spatial feature called ``central degree". This feature is easily recognized visually.

Let $P$ be the set of registration marks of an arrangement. A point $z$ in $P$ is called \emph{centrex} if every construction line that passes through $z$ is incident with the same number of registration marks on each ``side" of $z$. The number of construction lines through $z$ is called the \emph{central degree} of the arrangement. If a point $z$ does not exist, the arrangement does not have a central degree. In the following three arrangements, the centrex is highlighted in blue.

\begin{figure}[h!]
\centering
\includegraphics[width=4.6in]{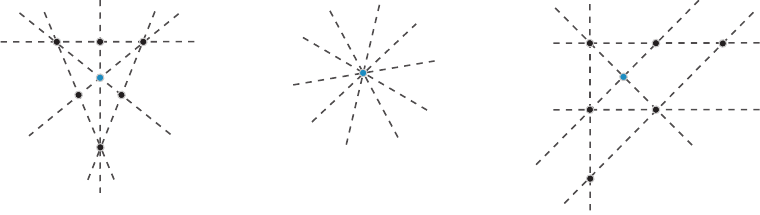}
\label{Figure20}
\end{figure}

A half-line originating at the centrex $z$ and passing through another registration mark of $P$ is called a \emph{spoke}. If the central degree of an arrangement is denoted by $d$, then there are 2$d$ spokes in total ``radiating" out of $z$. The spokes of an arrangement can be labelled in counter-clockwise order around the set $P$ as $C_1, ..., C_d, C_{d+1}, ..., C_{2d}$. The spokes $C_j$ and $C_{j + d}$ are opposite and have the same number $v_j = v_{j+d}$ of registration marks (the point $z$ is not included in this count). Then, the vector $(v_1, v_2, ..., v_d)$ is called the \emph{central signature} of the arrangement.

\vspace{12pt}

If two isomorphic arrangements have a centrex, then they can be distinguished in terms of the central signature defined by that centrex. If their central signatures are the same, then the arrangements cannot be distinguished on the basis of the comparison measures (1) through (4).

The triple of numbers in (1) is computed directly from the sets $L$ and $P$ of an arrangement. The measures in (2) and (3) are computed from the incidence structure that represents an arrangement \cite{haridis2020geometry}. The central degree and central signature in (4), however, cannot be computed from an incidence structure. Their computations depend on betweenness (i.e. order of collinear points), which is a property that incidence structures do not capture (hence the reason they can separate isomorphic arrangements successfully).

Finally, let us also note that affine and more generally projective transformations preserve collinearity and incidence. They carry an arrangement to another arrangement with the same incidence relation (this does not imply that the arrangements are actually equal, i.e. that their sets of registration marks and construction lines are equal). Thus, two arrangements that are projectively equivalent are also equivalent in terms of their incidence relation. Comparisons of arrangements by incidence is therefore ``coarser" than comparison by projective equivalence.

\end{document}